\documentclass[titlepage, 11pt]{article}
\usepackage{amssymb}
\usepackage{amscd} 
\usepackage{amsmath} 
\usepackage{amsthm}
\usepackage{color}
\usepackage{enumitem}
\usepackage{graphicx}
\usepackage[utf8]{inputenc} 
\usepackage{latexsym}
\usepackage{mathtools}
\usepackage{subcaption}
\usepackage{tikz}
\usepackage{mleftright}
\usepackage{comment}
\mleftright
\swapnumbers
\newtheorem{theorem}{Theorem}[section]
\newtheorem{lemma}[theorem]{Lemma} 
\newtheorem{prop}[theorem]{Proposition}

\newtheorem{conjecture}[theorem]{Conjecture}

\numberwithin{equation}{section}
\def\cupcup{\cup\cdots\cup}
\newcommand{\vare}{\varepsilon}
\newcommand\blackslug{\hbox{\hskip 1pt \vrule width 4pt height 8pt depth 1.5pt
        \hskip 1pt}}
\newcommand\bbox{\hfill \quad \blackslug \bigbreak}

\title{Clique covers of $H$-free graphs}
\author{
Tung Nguyen\thanks{Supported by AFOSR grant
A9550-19-1-0187.}\\
Princeton University, Princeton, NJ 08544, USA
\\
\\
Alex Scott\thanks{Research supported by EPSRC grant EP/X013642/1.}\\
Mathematical Institute, University of Oxford, Oxford OX2 6GG, UK
\\
\\
Paul Seymour\thanks{Supported by AFOSR grants
A9550-19-1-0187 and FA9550-22-1-0234, and NSF grant DMS-2154169.}\\
Princeton University, Princeton, NJ 08544, USA
\\
\\
St\'ephan Thomass\'e\\
Univ. Lyon, CNRS, ENS de Lyon, Universit\'e Claude Bernard Lyon 1, LIP UMR5668, France
}

\date{}
\newcommand{\Proof}{\noindent{\bf Proof.}\ \ }

\begin{document}

\maketitle
\begin{abstract}
It takes $n^2/4$ cliques to cover all the edges of a complete bipartite graph $K_{n/2,n/2}$, but how many cliques
does it take to cover all the edges of a graph $G$ if $G$ has no $K_{t,t}$ induced subgraph? We prove that 
$O(|G|^{2-1/(2t)})$ cliques suffice; and also prove that, even for graphs with no stable set of size four,
we may need more
than linearly many cliques. This settles two questions discussed at a recent conference in Lyon.
\end{abstract}

\section{Introduction}

A clique $X$ of a graph $G$ {\em covers} an edge $uv$ of $G$ if $u,v\in X$; and 
a {\em clique cover} of $G$ is a collection of cliques of $G$ that together cover all the edges.  
The {\em size} of a clique cover is the number of cliques in the collection.  What can we say about the sizes of clique covers?

The complete bipartite graph $K_{\lceil n/2\rceil,\lfloor n/2\rfloor}$ shows that, for an $n$-vertex graph, we may need as many as
$\lfloor n^2/4\rfloor$ cliques for a clique cover. In fact every graph $G$ has a clique cover of size at most $\lfloor |G|^2/4\rfloor$ 
(to 
see this, note that if $x,y$ are adjacent, we can cover all edges incident with $x$ or $y$ with 
at most $|G|-1$ cliques, so we may delete $x,y$ and use induction on $|G|$). But what if we restrict to $H$-free graphs? 
(A graph is {\em $H$-free} if it does not contain an induced copy of $H$.)
To make a difference, $H$ must be complete bipartite, or else $K_{\lceil n/2\rceil,\lfloor n/2\rfloor}$ is $H$-free; but what happens 
when $H$ is complete bipartite?

Indeed, what happens if $H=K_{s,0}$? Thus a graph $G$ is $H$-free if and only if $\alpha(G)<s$. (The sizes of the largest stable set 
and the
largest clique in $G$ are denoted by $\alpha(G), \omega(G)$ respectively.) The minimum size of clique covers in graphs $G$ with $\alpha(G)$
bounded already involves interesting questions.
For example, there is a long-standing conjecture that:

\begin{conjecture}\label{alpha2}
If $\alpha(G)\le 2$, there is a clique cover of size at most $|G|$.
\end{conjecture}
(``Long-standing'', but we do not know the source. Seymour recalls working on it many years ago, possibly in the 1980's.)
Which other graphs $H$ have the property that every $H$-free graph has a clique cover of size at most $|G|$?
It turns out that $H$ must be an induced subgraph of $K_{1,3}$.
(To see this, observe that every such graph $H$ must be an induced subgraph of a complete bipartite graph, and of 
the graph obtained from $K_{2,3}$ by subdividing two disjoint edges.) This leads us to the case when $H=K_{1,3}$, and for that there is 
a remarkable result of Javadi and Hajebi~\cite{hajebi}:
\begin{theorem}\label{clawfree}
If $G$ is connected and $K_{1,3}$-free, and has a stable set of cardinality three, then $G$ admits a clique cover of size at most $|G|$.
\end{theorem}
Thus, if \ref{alpha2} is true, then every $H$-free graph has a clique cover of size at most $|G|$ if and only if $H$ is an induced subgraph of $K_{1,3}$.

We have nothing to contribute to \ref{alpha2} itself, but what if we increase the bound on $\alpha(G)$?
Javadi and Hajebi~\cite{hajebi} asked whether 
all graphs $G$ with $\alpha(G)$ at most a constant admit clique covers of size $O(|G|)$, but we will disprove this.
We will show that:

\begin{theorem}\label{counterex}
There exists $C>0$ such that for infinitely many $n$, there is a graph $G$ on $n$ vertices with $\alpha(G)\le 3$
that
requires $Cn^{6/5}/(\log n)^2$ cliques in any clique cover.
\end{theorem}

And as an upper bound, we will show:
\begin{theorem}\label{stableup}
For every integer $s\ge 3$, if $G$ is a graph with no stable set of size $s$, then $G$ admits a clique cover of size at most
$O(|G|^{2-\frac{1}{s-1}})$.
\end{theorem}

At the other extreme, what happens if we exclude $K_{t,t}$? Sepehr Hajebi~\cite{lyon} recently proposed the following:
\begin{conjecture}\label{hajebiconj1}
For every integer $t\ge1$ there exists $\vare>0$ such that every $K_{t,t}$-free  $G$ has a clique cover of size $O(|G|^{2-\vare})$.
\end{conjecture}
Our main result is a proof of \ref{hajebiconj1}. We will show that:

\begin{theorem}\label{mainthm}\text{ }
\begin{itemize}
\item For all integers $s,t$ with $s\ge 3$ and $t\ge 2$, every $K_{s,t}$-free graph $G$ with sufficiently many vertices
has a clique cover of size at most $\frac32 |G|^{2-1/(s+t)}$.
\item For all integers $s$ with $s\ge 3$, 
every $K_{s,1}$-free graph $G$ (and {\rm a fortiori}, every  $K_{s,0}$-free graph)
has a clique cover of size at most 
$$O\left( \left(\frac{|G|}{\log |G|}\right)^{\frac{s-2}{s-1}} |G|\right).$$
\item Every $K_{2,2}$-free graph $G$ has a clique cover of size $O(|G|^{3/2})$;
\item Every $K_{2,3}$-free graph $G$ has a clique cover of size $O(|G|^{3/2}(\log |G|)^{1/2})$.
\end{itemize}
\end{theorem}
The second bullet implies \ref{stableup}.
We observe that the third bullet here is asymptotically sharp, since there are bipartite $K_{2,2}$-free graphs $G$ with 
$\Omega(|G|^{3/2})$ edges.

\section{Subquadratic clique covers}

In this section we will prove \ref{mainthm}. We begin with some lemmas. 
Ajtai,
Koml\'os and Szemer\'edi~\cite{ajtai} showed (logarithms in this paper are to base two):
\begin{lemma}\label{ajtai1}
For every integer $s\ge 2$ there exists $c>0$ such that, for all integers $a\ge 2$, the Ramsey number
$$R(a,s)\le \frac{ca^{s-1}}{(\log a)^{s-2}};$$ 
that is, every graph with at least $\frac{ca^{s-1}}{(\log a)^{s-2}}$ vertices has either a clique of size $a$ or a stable set of size $s$.
\end{lemma}
Let us rewrite \ref{ajtai1} in a form more convenient for us:
\begin{lemma}\label{ajtai2}
For every integer $s\ge 2$ there exists $c>0$ such that
if $w> 1$ is some real number, and $G$ is a graph with $\alpha(G)<s$ and $\omega(G)\le w$, then 
$$|G|< \frac{cw^{s-1}}{(\log w)^{s-2}}.$$
\end{lemma}
\Proof
Choose $c'$ such that setting $c=c'$ satisfies \ref{ajtai1}. Let $c=2^{s-1}c'$; we claim that $c$ satisfies \ref{ajtai2}.
Let $w\ge 1$, and let $G$ be a graph with $\alpha(G)<s$ and $\omega(G)\le w$. Let $a=\lfloor w\rfloor +1$. Then $a\ge 2$ is an integer, and
$\omega(G)<a$. By \ref{ajtai1}, 
$$|G|<\frac{c'a^{s-1}}{(\log a)^{s-2}}.$$
But $a\le 2w$ since $w\ge 1$, and so $c'a^{s-1}\le cw^{s-1}$, and since $(\log a)^{s-2}\ge (\log w)^{s-2}$, this proves \ref{ajtai2}.~\bbox

A theorem of Erd\H{o}s and Hajnal~\cite{EH}, in support of their well-known conjecture, implies that for all $s,t$ there exists $c>0$ such that if $G$ is $K_{s,t}$-free then $G$ has a clique or stable set 
of cardinality at least $|G|^c$. But we  want to make the result as sharp as we can, so we will do it again.

\begin{lemma}\label{EHt}
Let $s,t$ be integers with $s\ge t\ge 2$, and let $c\ge s$ satisfy \ref{ajtai2}.
If $G$ is $K_{s,t}$-free and $w>1$ is a real number with $w\ge \omega(G)$, then 
$$|G|\le \frac{c\alpha(G)^tw^{s-1}}{(\log w)^{s-2}}.$$
\end{lemma}
\Proof 
We may assume that $\alpha(G)\ge 2$, because otherwise 
$$|G|=\omega(G)\le \frac{c\alpha(G)^tw^{s-1}}{(\log w)^{s-2}}$$ 
(since $c\ge 1$ and $s\ge 2$, and $\omega(G)\le w$, and $w\ge \log w$), and the theorem holds.
\\
\\
(1) {\em $V(G)$ is the union of at most $\alpha(G)^t$ sets  each including no stable set of cardinality $s$.}
\\
\\
If $\alpha(G)<s$, the claim holds, so we may assume that $\alpha(G)\ge s$.
Let $S$ be a stable set of cardinality $\alpha(G)\ge s$. For $i\in \{t-1,t\}$, let $\mathcal{A}_{i}$ be the set of all subsets of $S$ of cardinality $i$.
For each $X\in \mathcal{A}_{t-1}$,
let $R_X$ be the set of all $v\in V(G)$ such that all neighbours of $v$ in $S$ belong to $X$ (thus, $X\subseteq R_X$). Since
$S$ is a largest stable set of $G$, it follows that $\alpha(G[R_X])\le t-1\le s-1$, because if there were a larger stable set in 
$R_X$, its union with $S\setminus X$ would be a stable set larger than $S$. 
For each $X\in \mathcal{A}_t$, let $R_X$ be the set of all $v\in V(G)\setminus S$ that are adjacent to every vertex in $X$.
Then $\alpha(G[R_X])\le s-1$ since $G$ is $K_{s,t}$-free. But every vertex with at most $t-1$ neighbours in $S$ belongs to $R_X$
for some $X\in \mathcal{A}_{t-1}$ (here we use that $|S|\ge t-1\ge 1$), and every vertex with at least $t$ neighbours in $S$ 
belongs to $R_X$ for some $X\in \mathcal{A}_t$,
and so $V(G)$ is the union of the sets $R_X\;(X\in \mathcal{A}_{t-1}\cup \mathcal{A}_t)$.
Moreover, 
$$|\mathcal{A}_{t-1}|+|\mathcal{A}_t| = \binom{\alpha(G)}{t-1}+ \binom{\alpha(G)}{t}=\binom{\alpha(G)+1}{t}\le \alpha(G)^t.$$
This proves (1).

\bigskip

From the choice of $c$,
if $R\subseteq V(G)$ includes no stable set of size $s$, then 
$$|R|\le \frac{cw^{s-1}}{(\log w)^{s-2}}.$$
By (1), it follows that 
$$|G|\le \frac{c\alpha(G)^tw^{s-1}}{(\log w)^{s-2}}.$$ 
This proves \ref{EHt}.~\bbox

We remark that the hypothesis $t\ge 2$ is not necessary. The same statement is true for $t=0,1$, but needs a slightly modified proof,
which we omit since we only need the result for $t\ge 2$.

\ref{EHt} implies that if $G$ is $K_{s,t}$-free then $\max(\alpha(G),\omega(G))\ge O\left(|G|^{1/(s+t-1)}\right)$.  
A similar proof shows that for every complete multipartite graph $H$, there exists $\vare$ such that if $G$ is $H$-free then 
$$\max(\alpha(G),\omega(G))\ge \vare |G|^{1/(|H|-1)}.$$ 
(We omit the proof, since we shall not use the result.)

Let us prove the first statement of \ref{mainthm}, the following:
\begin{theorem}\label{mainthm2}
Let $s\ge 3$ and $t\ge 2$ be integers. There exists $N$ such that every $K_{s,t}$-free graph $G$ with at least $N$ vertices 
admits a clique cover of size at most $\frac32 |G|^{2-1/(s+t)}$.
\end{theorem}
\Proof
We may assume that $s\ge t$, by exchanging them if necessary.
Let $c$ satisfy \ref{ajtai2}.
Since $s\ge 3$ we may choose $N$ such that
$$(\log N)^{s-2}\ge c(2t)^{t}(s+t)^{s-2}$$
(this is the only place in the proof that we need $s\ge 3$).
Let $d=1/(s+t)$, and let $G$ be a $K_{s,t}$-free graph with $n\ge N$ vertices.  We must show that $G$ has a clique cover of size at most $\frac32 n^{2-d}$;
and so we may assume that $\omega(G)\ge 2$.
We begin by choosing a maximal sequence of cliques 
in $G$, such that each clique covers at least $n^d$ edges not covered by previous cliques.  Thus, so far we have used at most 
$\frac12 n^{2-d}$ cliques.

Next, if there is any vertex $v$ that is incident with at most $n^{1-d}$ edges that have not yet been covered, we take copies 
of $K_2$ to cover all the uncovered edges incident with $v$.  Repeat this process until no such vertices remain.  Note that this 
step uses at most $n^{2-d}$ cliques in total, so altogether we have used at most $\frac32 n^{2-d}$ cliques. 

We claim that all edges of $G$ have now been covered; so, for a contradiction, suppose not.
Call a vertex $x$ {\em happy} if all edges incident with $x$ have been covered and {\em unhappy} otherwise;  thus there is at least one 
unhappy vertex.  Let $H$ be the subgraph of $G$ with vertex set the unhappy vertices and edge set the uncovered edges.  
Then $H$ has minimum degree at least $n^{1-d}$. Furthermore, no clique of $G$ covers at least $n^{d}$ edges of $H$, 
or we could have added it to our maximal sequence at the first step.  

Fix an unhappy vertex $v$, and let $D$ be the set of its neighbours in $H$, so $|D|\ge n^{1-d}$.  There is no clique
$K$ of $G[D]$ with size at least $n^d$, since adding $v$ to $K$ would give a clique of $G$ that covers $n^d$ edges of $H$ 
(all the edges from $v$ to $K$).  So by \ref{EHt}, taking $w=n^d$, it follows that
$G[D]$ contains a stable set $S$ where 
$$\frac{c|S|^t(n^d)^{s-1}}{(\log n^d )^{s-2}}\ge |D|\ge n^{1-d},$$
that is, 
$$|S|\ge d^{(s-2)/t}c^{-1/t} n^{(1-ds)/t}(\log n)^{(s-2)/t}=d^{(s-2)/t}c^{-1/t} n^{d}(\log n)^{(s-2)/t}$$
(since $d=(1-ds)/t$).

By a {\em copy} of $K_{1,t}$ we mean an induced subgraph of $G$ isomorphic to $K_{1,t}$; and a {\em leaf} of a graph
means a vertex with degree one.
We count copies of $K_{1,t}$ in $H$ with all their leaves in $S$. 
Let $L\subseteq S$ with $|L|=t$, and let $M$ be the set of vertices in $V(H)\setminus S$ that are adjacent in $H$ to every 
vertex in $L$. Since $L$ is stable in $G$ (as it is a subset of $S$), and $G$ is $K_{s,t}$-free, it follows that $M$ 
does not contain a stable set (of $G$) of size $s$. Moreover, $M$ contains no clique (of $G$) of size at least $n^d$, 
since adding any vertex of 
$L$ to such a clique would give a clique in $G$ covering at least $n^d$ edges from $H$. From \ref{ajtai2},  
$|M|\le \frac{c(n^d)^{s-1}}{(\log n^d)^{s-2}}$.
Since this holds for each choice of $L$, and there are only $\binom{|S|}{t}$ choices of $L$, it
follows that there are at most 
$$\frac{c(n^d)^{s-1}}{(\log n^d)^{s-2}}\binom{|S|}{t}\le \frac{c(n^d)^{s-1}}{(\log n^d)^{s-2}}\frac{|S|^t}{t!}$$
copies of $K_{1,t}$ with all leaves in $S$.

On the other hand, there are at least $n^{1-d}|S|$ edges of $H$ with an end in $S$.
For each $y\in V(H)$, let $r(y)$ be the set of vertices in $S$ adjacent in $H$ to $y$, and let $r$ be the average of the $r(y)$ over $y\in V(H)$. Thus
$$r \ge n^{1-d}|S|/|H|\ge n^{-d}|S|\ge d^{(s-2)/t}c^{-1/t}(\log n)^{(s-2)/t}.$$
Moreover, since $n\ge N$, it follows that 
$$(\log n)^{s-2}\ge c(2t)^{t}d^{2-s}$$
and so 
$$r\ge  d^{(s-2)/t}c^{-1/t}(\log n)^{(s-2)/t}\ge 2t.$$
The number of copies of $K_{1,t}$ with all leaves in $S$ is at least
$$\sum_{y\in V(H)}\binom{r(y)}{t}$$
(taking $\binom ab=0$ when $a<b$); and hence at least 
$$|H|\binom rt\ge |H|\frac{(r-t)^t}{t!}\ge |H|\frac{(r/2)^t}{t!}$$
by convexity and since $r\ge 2t\ge t$.
Consequently 
$$ |H|\frac{(r/2)^t}{t!} \le \frac{c(n^d)^{s-1}}{(\log n^d)^{s-2}}\frac{|S|^t}{t!},$$
that is,
$$ |H|(r/2)^t \le \frac{c(n^d)^{s-1}}{(\log n^d)^{s-2}}|S|^t.$$
Since $|S|\le r|H|n^{d-1}$ and $d(s+t)=1$, the right side of the above
is at most 
$$\frac{c(n^d)^{s-1}}{(\log  n^d )^{s-2}}r^t|H|^tn^{t(d-1)}=\frac{c}{d^{s-2}(\log n)^{s-2}}r^t|H|^tn^{1-d-t}\le 
\frac{c}{d^{s-2}(\log n)^{s-2}}r^t|H|n^{-d}.$$
Consequently
$$|H|(r/2)^t\le \frac{c}{d^{s-2}(\log n)^{s-2}}r^t|H|n^{-d}.$$
that is,
$$(\log n)^{s-2}n^{d}d^{s-2}\le  c 2^t,$$
contradicting that $n\ge N$.
This proves \ref{mainthm2}.~\bbox

A simplified version of the same argument yields a weakened form of the second statement of \ref{mainthm}:
\begin{prop}\label{stableup2}
For every integer $s\ge 3$ and $t\in \{0,1\}$, if $G$ is a $K_{s,t}$-free graph, then $G$ admits a clique cover of size at most
$O(|G|^{2-1/s})$.
\end{prop}
\Proof
Let $c$ be as in the proof of \ref{mainthm2}, let $d=1/s$, and choose $N$ such that $(\log N)^{s-2}> cs^{s-2}$. We will show that if $|G|=n\ge N$ and $G$ is $K_{s,t}$-free then $G$ admits a clique cover of size at most
$\frac32 |G|^{2-1/s}$, which implies the result. As in the proof of \ref{mainthm2}, we may assume that $G$ has a nonnull subgraph $H$ with minimum degree at least 
$n^{1-d}$, such that no clique of $G$ covers $n^d$ edges of $H$. Choose $v,D$ as before; then $G[D]$ has no clique of size $n^d$,
and no stable set of size $s$ (because $G$ is $K_{s,t}$-free and $t\le 1$), and so by \ref{ajtai2}, 
$$|D|<\frac{cn^{d(s-1)}}{(d\log n)^{s-2}}.$$
But $|D|\ge n^{1-d}$, and so 
$$n^{1-d}\le \frac{cn^{d(s-1)}}{(d\log n)^{s-2}}.$$
Since $1-d=d(s-1)$, it follows that 
$(\log n)^{s-2}\le cs^{s-2}$,
contradicting that $n\ge N$. This proves \ref{stableup2}.~\bbox

But we can do a little better. To prove the second statement of \ref{mainthm} as stated, we use a consequence of \ref{ajtai2}:
\begin{lemma}\label{ajtai5}
For all integers $s\ge 2$ there exists $c>0$ such that
if $G$ is a graph with $\alpha(G)<s$ and $|G|\ge 2$, then
$$|G|< \frac{c\omega(G)^{s-1}}{(\log |G|)^{s-2}}.$$
\end{lemma}
\Proof
Choose $d$ such that \ref{ajtai2} holds with $c=d$. We may assume that $d\ge s^{1/2}$ by increasing $d$.
Choose $c$ such that $c\ge d^2(2\log d)^{s-2}$ and $c>d(2(s-1))^{s-2}$; we will show that $c$ satisfies the lemma.

Let $G$ be a graph with $\alpha(G)<s$ and $|G|\ge 2$. If $\log|G|\le 2\log d $, then $|G|\le d^2$, and so 
$$|G|(\log |G|)^{s-2}\le d^2(2\log d)^{s-2}\le c\le c\omega(G)^{s-1}$$
as required. Thus we may assume that $\log|G|> 2\log d $. In particular $|G|\ge d^2\ge s$ and so $G$ has an edge.
By \ref{ajtai2}, 
$$|G|< \frac{d\omega(G)^{s-1}}{(\log \omega(G))^{s-2}}\le d\omega(G)^{s-1}$$
and so $\log|G|\le \log(d) +(s-1) \log(\omega(G))$.
Since $\log|G|>2\log(d)$, it follows that $\log|G|<2(s-1) \log(\omega(G))$.
Hence
$$|G|< \frac{d\omega(G)^{s-1}}{(\log \omega(G))^{s-2}}\le \frac{d(2(s-1))^{s-2}\omega(G)^{s-1}}{(\log |G|)^{s-2}}\le \frac{c\omega(G)^{s-1}}{(\log |G|)^{s-2}}$$
as required. This proves \ref{ajtai5}.~\bbox

We deduce:

\begin{lemma}\label{ajtai6}
For all integers $s\ge 3$ there exists $c>0$ such that
for every graph $G$ with $\alpha(G)<s$ and $|G|\ge 2$, $V(G)$ is the union of at most
$$c\left(\frac{|G|}{\log |G|}\right)^{\frac{s-2}{s-1}}$$
cliques.
\end{lemma}
\Proof Choose $d$ such that \ref{ajtai5} holds (with $c$ replaced by $d$).
Choose $f$ such that 
$2d f^{s-1}\ge (1/2)^{s-2}$.
Choose $N\ge 4$ such that $\log N>(1-2^{-\frac{1}{s-2}})^{-1}$,
and choose $c$ such that $c\ge 2/f$, and $c(n/\log n)^{\frac{s-2}{s-1}}\ge n$ for all nonzero integers $n\le N$.
We will show that
$c$ satisfies \ref{ajtai6}. Let $G$ be a graph with $\alpha(G)<s$.
We prove that the statement of the theorem is true for $G$, by induction on $|G|$.
If $|G|\le N$, then $V(G)$ is the union of $|G|\le c(|G|/\log |G|)^{\frac{s-2}{s-1}}$
cliques and the theorem holds, so we may assume that $|G|>N$.

Choose as many pairwise disjoint cliques as possible that each have cardinality at least
$$f|G|^{\frac{1}{s-1}}(\log |G|)^{\frac{s-2}{s-1}},$$ 
say $A_1\ldots A_k$. 
Let $G'=G\setminus (A_1\cupcup A_k)$. Thus 
$$\omega(G')< f|G|^{\frac{1}{s-1}}(\log |G|)^{\frac{s-2}{s-1}}.$$
\noindent(1) {\em $|G'|\le |G|/2$.}
\\
\\
Suppose not; then by \ref{ajtai5}, 
$$|G'|< \frac{d\omega(G')^{s-1}}{(\log |G'|)^{s-2}},$$
and so 
$$(|G|/2)(\log (|G|/2))^{s-2}\le |G'|(\log |G'|)^{s-2}< d\omega(G')^{s-1}\le d f^{s-1}|G|(\log |G|)^{s-2}.$$
Thus 
$$(\log |G|-1))^{s-2}\le 2d f^{s-1}(\log |G|)^{s-2}.$$
But $\log |G| -1 \ge \frac12 \log |G|$
(because $|G|\ge N\ge 4$),
and so 
$(1/2)^{s-2}\le 2d f^{s-1}$,
a contradiction. This proves (1).

\bigskip

Since $A_1\ldots A_k$ all have cardinality at least $f|G|^{\frac{1}{s-1}}(\log |G|)^{\frac{s-2}{s-1}}$, it follows that
$$k\le \frac{f^{-1}|G|^{\frac{s-2}{s-1}}}{(\log |G|)^{\frac{s-2}{s-1}}}=f^{-1}\left(\frac{|G|}{\log |G|}\right)^{\frac{s-2}{s-1}}\le (c/2)\left(\frac{|G|}{\log |G|}\right)^{\frac{s-2}{s-1}}.$$

From the inductive hypothesis, if $|G'|\ge 2$ then $V(G')$ is the union of 
$$c\left(\frac{|G'|}{\log |G'|}\right)^{\frac{s-2}{s-1}}$$
cliques. Hence, $V(G')$ is the union of at most 
$$c\left(\frac{|G|/2}{\log |G|-1}\right)^{\frac{s-2}{s-1}}$$
cliques, by (1), even if $|G'|\le 1$. 
But
$$c\left(\frac{|G|/2}{\log |G|-1}\right)^{\frac{s-2}{s-1}}\le (c/2)\left(\frac{|G|}{\log |G|}\right)^{\frac{s-2}{s-1}}$$
since 
$$2^{\frac{1}{s-2}}\ge \frac{\log |G|}{\log |G|-1}.$$
Thus, both $A_1\cupcup A_k$ and $V(G')$ are the union of at most 
$$(c/2)\left(\frac{|G|}{\log |G|}\right)^{\frac{s-2}{s-1}}$$
cliques. 
Adding, this proves \ref{ajtai6}.~\bbox

We use this to show the second statement of \ref{mainthm}, the following:
\begin{theorem}\label{stableup4}
For every integer $s\ge 3$,  let $c$ be as in \ref{ajtai6}; if $G$ is a $K_{s,1}$-free graph with $|G|\ge 2$, 
then $G$ admits a clique cover of size at most
$$\frac{c|G|^{2-\frac{1}{s-1}}}{(\log |G|)^{\frac{s-2}{s-1}}}.$$
\end{theorem}
\Proof
By \ref{ajtai6}, there is a set $\mathcal{A}$ of cliques of $G$ with union $V(G)$
and with 
$$|\mathcal{A}|\le c\left(\frac{|G|}{\log |G|}\right)^{\frac{s-2}{s-1}}.$$
For each $v\in V(G)$ and $A\in \mathcal{A}$, let $A_v$ be the clique consisting of $v$ and the set of neighbours
of $v$ that belong to $A$. Then the set of all the cliques $A_v$ is a clique cover satisfying the theorem.
This proves \ref{stableup4}.~\bbox

Now we prove the third and fourth statements of \ref{mainthm}.
We will need the following, which is implied by \ref{ajtai2} with $s=3$:
\begin{lemma}\label{ajtai}
There exists $k>0$ such that every graph $G$ with no stable set of size three has a clique of size at least
$k|G|^{1/2}\sqrt{\log|G|}$.
\end{lemma}

We will show:
\begin{theorem}\label{Kst}\text{ }
\begin{itemize}
\item Every $K_{2,2}$-free graph $G$ has a clique cover of size $O(|G|^{3/2})$.
\item Every $K_{2,3}$-free graph $G$ has a clique cover of size $O(|G|^{3/2}(\log |G|)^{1/2})$.
\end{itemize}
\end{theorem}
\Proof
The proofs for both statements are much the same, and we will do them at the same time.
Let $G$ be either $K_{2,2}$-free or $K_{2,3}$-free, 
let $v$ be a vertex of minimum degree, and let $D$ be the set of its neighbours. 
We will show that $D$ is the union of a small
number of cliques.  Adding $v$ to each of these cliques, we see that the edges incident with $v$ 
can be covered by the same small number of cliques; thus we may delete $v$ and argue by induction.
It remains to show that $D$ is the union of an appropriately small
number of cliques.

First we need:
\\
\\
(1) {\em Let $M\subseteq D$.
Then either:
\begin{itemize}
\item there is a set $J_1\subseteq M$ with $|J_1|\ge |M|^2/(4n)$, and two nonadjacent vertices $x,y\in V(G)\setminus J_1$, both adjacent to 
every vertex in $J_1$; or
\item there is a clique $J_2\subseteq M$ with $|J_2|\ge |M|/4$.
\end{itemize}
}

\bigskip

\noindent Let $A\subseteq M$ be the set of vertices in $M$ with at least $|M|/3$ neighbours outside $D\cup\{v\}$, 
and let $B=M\setminus A$.
Suppose first that $|A|\ge 3|M|/4$. Then the number of edges from $A$ to $V(G)\setminus (D\cup\{v\})$ is at least $|M|^2/4$,
and so some vertex $x \in V(G)\setminus (D\cup\{v\})$ 
has a set $J$ of at least $|M|^2/(4n)$ neighbours in $M$, and the first bullet of (1) holds (taking $y=v$).

Otherwise $|B|\ge |M|/4$.  If $B$ is a clique then the second bullet holds. Otherwise there are nonadjacent vertices $x,y\in B$;
and as $x,y$ each have at most
$|M|/3$ non-neighbours in $D$ (because $v$ was chosen with minimum degree, and $x,y\in B$), there are at most $2|M|/3$
vertices in $M$ nonadjacent to one of $x,y$ (counting $x,y$ themselves);
and so $x,y$ have at least $|M|/3$ common neighbours
in $M$, and the first bullet holds. This proves (1).

\bigskip

We deduce:
\\
\\
(2) {\em If $G$ is $K_{2,2}$-free, then for every $M\subseteq D$, there is a clique in $M$ with size at least 
$|M|^2/(4n)$.}
\\  
\\  
This is immediate from (1), because the set $J_1$ in (1) must be a clique, since $G$ is $K_{2,2}$-free, and
the set $J_2$ satisfies $|J_2|\ge |M|/4 \ge |M^2/(4n)$. This proves (2).
\\
\\
(3) {\em Let $k$ satisfy \ref{ajtai}, and let    
$\beta=\min(k/2, 1/4)$. If $G$ is $K_{2,3}$-free,
then for every $M\subseteq D$ with $|M|^2\ge 4n$, there is a clique in $M$ with size at least
$$\frac{\beta|M|}{\sqrt{n}}\sqrt{\log (|M|^2/(4n))}.$$}

\noindent By (1), one of the sets $J_1,J_2$ of (1) exist. If $J_1$ exists, then it contains no stable triple of vertices, and so
by \ref{ajtai}, it contains a clique of size at least 
$$k (|J_1|\log(|J_1|))^{1/2}\ge  \frac{k |M|}{2\sqrt n } (\log (|M|^2/(4n)))^{1/2},$$
and the claim holds. 
If $J_2$ exists then again the claim holds since 
$$|M|/4\ge \beta|M|\ge \frac{\beta|M|}{\sqrt{n}}\sqrt{\log (|M|^2/(4n))}.$$
This proves (3).

\bigskip
Now we will use (2) or (3) to show 
that the vertices in $D$ can be covered by an
appropriately small collection $\mathcal C$ of cliques.  
We choose $\mathcal C$ by choosing greedily a largest clique among the uncovered vertices of $D$ until at most $4\sqrt n$ vertices 
remain,  and then covering the remaining vertices by singletons. To bound the total number of cliques, we track the process, writing 
$M$ for the set of uncovered vertices at each stage. We divide the values of $|M|$ into ranges $[1,4\sqrt n)$ and 
$[2^i\sqrt n,2^{i+1}\sqrt n)$ for $i\ge 2$. 

We assume first that $G$ is $K_{2,2}$-free.
Thus by (2), if $|M|$ is in the range  $[2^i\sqrt n,2^{i+1}\sqrt n)$, then the size of the clique we obtain is at least 
$|M|^2/(4n) \ge  2^{2i-2}$, and so there will
be at most 
$$\frac{2^{i+1}\sqrt{n}}{2^{2i-2}}= \frac{8\sqrt{n}}{2^i}$$ 
cliques chosen for $|M|$ in this range.    
The total number of cliques in $\mathcal C$ is therefore at most
$$
\sum_{i\ge 2}\frac{8\sqrt{n}}{2^{i}} + 4\sqrt n =O(\sqrt n).
$$
Consequently the first bullet of the theorem follows by induction.

Now we assume that $G$ is $K_{2,3}$-free, and use (3) in place of (2). 
If $|M|$ is in the range  $[2^i\sqrt n,2^{i+1}\sqrt n)$ where $i\ge 2$, then the size of the clique we obtain is at least           
$$\frac{\beta|M|}{\sqrt n}\sqrt{\log (|M|^2/(4n))}\ge \beta 2^i \sqrt{\log(2^{2i-2})}=\beta 2^i \sqrt{2i-2}.$$
Consequently, at most 
$$\frac {2^{i+1}\sqrt n}{\beta 2^i \sqrt{2i-2}}=\frac {2\sqrt n}{\beta \sqrt{2i-2}}$$
cliques will be chosen during this range. 
Thus the total number of cliques is at most
$$
\sum_{i= 2}^{\log n}\frac{2\sqrt{n}}{\beta\sqrt{2i-2}} + 4\sqrt n =O(\sqrt {n\log n}).
$$
Hence the second bullet of the theorem follows by induction. This proves \ref{Kst}.~\bbox

\section{Lower bounds}

What can we say from the other side? For $K_{s,0}$-free graphs, the result of this section, with \ref{mainthm}, shows that 
(roughly speaking) the answer is somewhere between
$n^{2-4/(s+1)} $ and $n^{2-\frac{1}{s-1}}$. We need the following result of Spencer (theorem 2.2 of~\cite{spencer}):
\begin{lemma}\label{spencer}
For all integers $s\ge 3$, there exists $c>0$ such that for all integers $t\ge 3$, the Ramsey number
$R(s,t)$ is at least 
$c(t/\log t)^{\frac{s+1}{2}}.$
Consequently,  for all $s\ge 3$ there exists $C>0$ such that for infinitely many $n$,
there is a graph $J$ with $n$ vertices such that $\omega(G)< Cn^{\frac{2}{s+1}}\log n$ and $\alpha(G)<s$.
\end{lemma}

\begin{theorem}\label{betterex}
For all $s\ge 3$, there exists $c>0$ such that for infinitely many $n$, there is a graph with $n$ vertices and with no stable set of 
size $s$,
such that every clique cover has size at least $cn^{2-4/(s+1)}/(\log n)^2$.
\end{theorem}
\Proof
Choose $C$ as in the second statement of \ref{spencer}, and let $c$ satisfy $c^{-1}=C^22^{2-4/(s+1)}$. Now choose $m>0$ such that 
there is a graph $J$ with $m$ vertices, and with $\omega(J)< Cn^{\frac{2}{s+1}}\log m$ and $\alpha(J)<s$.
Let $n=2m$. 
Take two vertex-disjoint copies $J_1,J_2$ of $J$, and make every vertex of $J_1$
adjacent to every vertex of $J_2$, forming $G$; thus $|G|=n$. Then $G$ has no stable set of size $s$; 
and every clique of $G$ covers at most 
$C^2m^{\frac{4}{s+1}}(\log m)^2$ of the edges between $V(J_1)$ and $V(J_2)$. Since there are $m^2$ such edges, every clique cover of $G$
has size at least 
$$C^{-2}m^{2-\frac{4}{s+1}}/(\log m)^2\ge cn^{2-\frac{4}{s+1}}/(\log n)^2.$$
This proves \ref{betterex}.~\bbox

Taking $s=4$, this proves \ref{counterex}.

\end{document}